\documentclass[a4paper,12pt]{amsart}
\usepackage{amssymb}
\usepackage{amsmath}
\usepackage{ifthen}
\usepackage{hyperref}
\allowdisplaybreaks
\newtheorem{defn}{Definition}[section]
\newtheorem{thm}[equation]{Theorem}
\newtheorem{lem}[equation]{Lemma}

\begin{document}
\newcommand{\A}{{\mathcal A}}
\newcommand{\es}{{\mathcal S}}
\newcommand{\IR}{{\mathbb R}}
\newcommand{\IC}{{\mathbb C}}
\newcommand{\IN}{{\mathbb N}}
\newcommand{\D}{{\mathbb D}}
\newcommand{\M}{{\mathcal M}}
\newcommand{\N}{{\mathcal N}}
\newcommand{\CC}{{\mathcal C}}
\newcommand{\R}{{\mathcal R}}
\newcommand{\V}{{\mathcal V}}

\newcommand{\beq}{\begin{eqnarray}}
\newcommand{\beqq}{\begin{eqnarray*}}
\newcommand{\eeq}{\end{eqnarray}}
\newcommand{\eeqq}{\end{eqnarray*}}
\newcommand{\bdefe}{\begin{defn}}
\newcommand{\edefe}{\end{defn}}
\newcommand{\bthm}{\begin{thm}}
\newcommand{\ethm}{\end{thm}}
\newcommand{\blem}{\begin{lem}}
\newcommand{\elem}{\end{lem}}

\title[Univalent Functions with Non-Negative Coefficients involving Clausen's Hypergeometric Function]{Univalent Functions with Non-Negative Coefficients involving Clausen's Hypergeometric Function}

\author[K. Chandrasekran]{K. Chandrasekran}
\address{K. Chandrasekran \\ Research Scholar \\ Department of Mathematics \\ MIT Campus, Anna University \\ Chennai 600 044, India}
\email{kchandru2014@gmail.com}

\author[G. Murugusundaramoorthy]{G. Murugusundaramoorthy}
\address{G. Murugusundaramoorthy \\ School of Advanced Sciences\\ Vellore Institute of Technology\\ Vellore-632014, India}
\email{gmsmoorthy@yahoo.com}

\author[D. J. Prabhakaran ]{D. J. Prabhakaran}
\address{D. J. Prabhakaran \\ Department of Mathematics \\ MIT Campus, Anna University \\ Chennai 600 044, India}
\email{asirprabha@gmail.com}


\begin{abstract}
In this work, we derived the necessary and sufficient conditions on parameters for $_3F_2(^{a,b,c}_{b+1,c+1};z)$ Hypergeometric Function to be in the classes $\M^{\ast}(\lambda,\alpha)$ and $\N^{\ast}(\lambda,\alpha)$ and information regarding the image of function $_3F_2(^{a,b,c}_{b+1,c+1};z)$ belonging to $\R^{\tau}(A,B)$ by applying the convolution operator in open unit disc $\D =\{z:\, |z|<1\}$.
\end{abstract}
\maketitle

\maketitle
\section{Introduction}
Let $\D = \{z \in \IC: |z|<1\}$ be the open unit disc in the complex plane $\IC$. Let $\mathcal{H}$ denote the class of all analytic functions in $\D$. Let $\A$ denote the family of analytic functions $f$ of the form
\beq \label{chp2inteq0}
f(z)= z+\sum_{n=2}^{\infty}\, a_n\,z^n,\, z \in \D
\eeq with $f(0)=0$ and $f^{\prime}(0)=1$ in the open unit disc $\D$. Which is the subclass of $\mathcal{H}$ and Let, $ {\es} \subset \A$, \, i.e. $\es$ denotes the class of all normalised functions that are analytic and univalent in open unit disc $\D.$  For the function $\displaystyle f$ is given by (\ref{chp2inteq0}) in $\A$ and $g \in \A$  with $ \displaystyle g(z)= z+\sum_{n=2}^{\infty}\, b_n\,z^n $, the \emph{convolution product} of $f$ and $g$ is defined by $$\displaystyle (f*g)(z)= z+\sum_{n=2}^{\infty}\, a_n\,b_n\, z^n, z \in \D.$$
Note that the \emph{convolution product} is called \emph{Hadamard Product}. For more details refer \cite{chp2A-W-Goodman-1983-book}
\begin{defn}
The subclass $\V$ of $\A$ consisting of functions of the form
\begin{eqnarray*}
  f(z)= z+\sum_{n=2}^{\infty}\, a_n\,z^n, \, z \in \D,\,\, {\rm with}\, \, a_n \geq 0,\, n\in \IN,\,n \geq 2.
\end{eqnarray*}
\end{defn}

In \cite{chp2Uralegaddi-1994}, Uralegaddi et al. introduced the following two classes which are stated as:
\begin{defn}\cite{chp2Uralegaddi-1994}
The class $\M(\alpha)$ of \emph{starlike functions of order $\alpha$,} with $1<\alpha \leq \frac{4}{3}$, defined by
$$\M(\alpha) = \left\{f\in \mathcal{A}: \Re\left(\frac{zf^{\prime}(z)}{f(z)}\right)< \alpha,\, z \in \mathbb{D}\right\}$$
\end{defn}
\begin{defn}\cite{chp2Uralegaddi-1994}
The class $\N(\alpha)$ of \emph{convex functions of order $\alpha$,} with $1<\alpha \leq \frac{4}{3}$, defined by
$$\N(\alpha) = \left\{f\in \mathcal{A}: \Re\left(1+\frac{zf^{\prime\prime}(z)}{f^{\prime}(z)}\right)< \alpha,\, z \in \mathbb{D}\right\}=\
\left\{f \in \A : zf^{\prime}(z)\in \M{(\alpha)}\right\}$$
\end{defn}

In this paper, we considere the two subclasses $\M(\lambda, \alpha)$ and $\N(\lambda, \alpha)$ of  to discuss some inclusion properties based on Clausen's Hypergeometric Function. These two subclasses was introduced by Bulboaca and Murugusundaramoorthy \cite{chp2Bulboaca-Murugu-2020}. which are stated as follows:

\begin{defn}\cite{chp2Bulboaca-Murugu-2020} For some $\alpha\, \left(1<\alpha\leq \frac{4}{3}\right)$ and $\lambda\, \left(0\leq \lambda <1\right)$, the functions of the form (\ref{chp2inteq0}) be in the subclass $\M(\lambda,\alpha)$ of $\es$ is
  \begin{eqnarray*}
    \M(\lambda,\alpha) &=& \left\{ f \in \A:\Re\left(\frac{zf^{\prime}(z)}{(1-\lambda)f(z)+\lambda z\, f^{\prime}(z)}\right)< \alpha,\, z \in \D\right\}
  \end{eqnarray*}
\end{defn}
\begin{defn}\cite{chp2Bulboaca-Murugu-2020} For some $\alpha\, \left(1<\alpha\leq \frac{4}{3}\right)$ and $\lambda\, \left(0\leq \lambda <1\right)$, the functions of the form (\ref{chp2inteq0}) be in the subclass $\N(\lambda,\alpha)$ of $\es$ is
\begin{eqnarray*}
 \N(\lambda,\alpha) &=& \left\{f \in \A: \Re\left(\frac{f^{\prime}(z)+zf^{\prime\prime}(z)}{f^{\prime}(z)+\lambda z\, f^{\prime\prime}(z)}\right)< \alpha,\, z \in \D \right\}
\end{eqnarray*}
\end{defn}

Also, let $\M^{\ast}(\lambda,\alpha)\equiv \M(\lambda,\alpha)\cap \V$ and $\N^{\ast}(\lambda,\alpha)\equiv \N(\lambda,\alpha)\cap \V$.

\begin{defn}\label{dp01}\cite{chp2dixitpal1995}
  A function $f\in \A$ is said to be in the class $\R^{\tau}(A,B)$, with $\tau\in \IC\backslash\{0\}$ and $-1\leq B \leq A\leq 1$, if it satisfies the inequality $$\displaystyle\biggl|\frac{f^{\prime}(z)-1}{(A-B)\tau-B[f^{\prime}(z)-1]}\biggr|<1, z \in \D$$
\end{defn}

Dixit and Pal \cite{chp2dixitpal1995} introduced the Class $\R^{\tau}(A,B)$. Which is stated as in the definition \ref{dp01}. If we substitute $\tau=1,\, A=\beta\,$ and $ B=-\beta, \, (0 < \beta \leq 1)$ in the definition \ref{dp01},  then we obtain the class of functions $f \in \A$ satisfying the inequality
$$\biggl|\frac{f^{\prime}(z)-1}{f^{\prime}(z)+1}\biggr|<\beta,\, z \in \D$$
which was studied by Padmanabhan \cite{chp2padma1970} and others subsequently.

\begin{defn}\cite{chp2Andrews-Askey-Roy-1999-book}
The $_3F_2(a,b,c;d,e;z)$ hypergeometric series is defined as
\beq\label{chp2inteq5}
_3F_2(a,b,c;d,e;z)=\sum_{n=0}^{\infty}\frac{(a)_n(b)_n(c)_n}{(d)_n(e)_n(1)_n}z^n,\, \, \, a,b,c,d,e\in \IC,
\eeq provided $d,\, e\, \neq 0,-1,-2,-3\cdots,$ which is an analytic function in open unit disc $\D$.
\end{defn}
We consider the linear operator $\mathcal{I}^{a,b,c}_{b+1,c+1}(f):\A \rightarrow \A$ defined by convolution product
\beq\label{chp2inteq7}
\mathcal{I}^{a,b,c}_{b+1,c+1}(f)(z) = z\, _3F_2(^{a,b,c}_{b+1,c+1};z)*f(z) = z+\sum_{n=2}^{\infty} A_n\, z^n
\eeq
where $A_1=1$ and for $n > 1,$
\beq\label{chp2inteq007}
A_n&=&\frac{(a)_{n-1}(b)_{n-1}(c)_{n-1}}{(b+1)_{n-1}(c+1)_{n-1}(1)_{n-1}}\, a_n.
\eeq

Motivated by the results in connections between various subclasses of analytic univalent functions, by using hypergeometric functions   \cite{chp2chandru-prabha-Thesis-2022,chp2Chandru-prabha-2019,chp2Chandru-prabha-2020,chp2Chandru-prabha-2021,chp2Chandru-prabha-2022,chp2Uralegaddi-1994},  and Poisson distributions \cite{chp2Bulboaca-Murugu-2020}, we obtain the necessary and sufficient conditions on parameters for $_3F_2(^{a,b,c}_{b+1,c+1};z)$ hypergeometric series to be in the classes $\M^{\ast}(\lambda,\alpha)$ and $\N^{\ast}(\lambda,\alpha)$ and information regarding the image of functions $_3F_2(^{a,b,c}_{b+1,c+1};z)$ hypergeometric series belonging to $\R^{\tau}(A,B)$ by applying the \emph{Hadamard product.}

\section{Main Results and Proofs}
First, we recall the following results to prove our main theorems.
\blem\label{chp2lem1eqn1} \cite{chp2Murugu-2018}
For some $\alpha\, (1<\alpha \leq \frac{4}{3})$ and $\lambda \, (0\leq \lambda < 1)$, and if $f \in \V$, then $f \in \M^{\ast}(\lambda,\alpha)$ if and only if
\begin{eqnarray}\label{chp2lem1eqn2}
  \sum_{n=2}^{\infty}\, [n-(1+n\lambda-\lambda)\alpha]a_n &\leq& \alpha-1.
\end{eqnarray}
\elem
\blem\label{chp2lem2eqn1}\cite{chp2Murugu-2018}
For some $\alpha\, (1<\alpha \leq \frac{4}{3})$ and $\lambda \, (0\leq \lambda < 1)$, and if $f \in \V$, then $f \in \N^{\ast}(\lambda,\alpha)$ if and only if
\begin{eqnarray}\label{chp2lem2eqn2}
  \sum_{n=2}^{\infty}\, n\,[n-(1+n\lambda-\lambda)\alpha]a_n &\leq& \alpha-1.
\end{eqnarray}
\elem
The following result is due to  Miller and  Paris \cite{chp2Mill-Paris-2012-ITSF} $\&$ Shpot and Srivastava \cite{chp2Shpot-Srivas-2015-AMC}.
\begin{thm}
For $a,\, b,\,c > 0,\, c\neq b$ and $a < \min(1,\, b+1,\,c+1)$,
\beq\label{chp2inteq6}
 _3F_2\left(^{a,b,c}_{b+1,c+1};1\right)&=& \frac{bc}{c-b}\Gamma(1-a)\left[\frac{\Gamma(b)}{\Gamma(1-a+b)}-\frac{\Gamma(c)}{\Gamma(1-a+c)}\right].
\eeq
\end{thm}
Now, we state the following lemma due to Chandrasekran and Prabhakaran \cite{chp2Chandru-prabha-2019} which is useful to prove our main results.
\blem \label{chp2lem3eqn1} \cite{chp2Chandru-prabha-2019}
Let $a,b,c > 0$. Then we have the following:
\begin{enumerate}
\item For $ b, c > a-1$, we have
\begin{eqnarray*}
\sum_{n=0}^{\infty} \frac{(n+1)(a)_n\, (b)_n\, (c)_n }
{(b+1)_n\, (c+1)_n\, (1)_n} &=&  \frac{bc\, \Gamma(1-a) }{c-b}\left[\frac{(1-b)\Gamma(b)}{\Gamma(1-a+b)}-\frac{(1-c)\Gamma(c)}{\Gamma(1-a+c)}\right].
\end{eqnarray*}
\item For $ b, c > a-1$, we have
\begin{eqnarray*}
\sum_{n=0}^{\infty} \frac{(n+1)^2(a)_n\, (b)_n\, (c)_n }
{(b+1)_n\, (c+1)_n\, (1)_n} &=&  \frac{bc\, \Gamma(1-a)}{c-b}\left[\frac{(1-b)^2\Gamma(b)}{\Gamma(1-a+b)}-\frac{(1-c)^2\Gamma(c)}{\Gamma(1-a+c)}\right].
\end{eqnarray*}
\item For $ b, c > a-1$, we have
\begin{eqnarray*}
\sum_{n=0}^{\infty} \frac{(n+1)^3(a)_n\, (b)_n\, (c)_n }
{(b+1)_n\, (c+1)_n\, (1)_n} &=&  \frac{bc\, \Gamma(1-a)}{c-b} \left[\frac{(1-b)^3\Gamma(b)}{\Gamma(1-a+b)}-\frac{(1-c)^3\Gamma(c)}{\Gamma(1-a+c)}\right].
\end{eqnarray*}
\item For $a\neq 1,\, b\neq 1,\,$ and $c\neq 1$ with $\, b, c >\max\{0,  a-1\}$, we have
\begin{eqnarray*}
\sum_{n=0}^{\infty} \frac{(a)_n\, (b)_n\, (c)_n }
{(b+1)_n\, (c+1)_n\, (1)_{(n+1)}} &=&  \frac{bc}{(a-1)(b-1)(c-1)}\\ && \times\left[\frac{\Gamma(2-a)}{c-b}\left(\frac{(c-1)\Gamma(b)}{\Gamma(1-a+b)}-\frac{(b-1)\Gamma(c)}{\Gamma(1-a+c)}\right)-1\right].
\end{eqnarray*}
\end{enumerate}
\elem

\bthm \label{chp2thm1eqn1}
 Let $a \in {\IC} \backslash \{ 0 \}, \, b,\,c > 0$,\, $c\neq b$  and $ |a| < \min\{1,\, b+1, \, c+1\}$. A sufficient condition for the function $z\, _3F_2\left(^{a,b,c}_{b+1,c+1};z\right) $ to belong to the class $\M^{\ast}(\lambda,\alpha),\, 1<\alpha \leq \frac{4}{3}$ and $ 0 \leq \lambda  < 1 $ is that
\beq\label{chp2thm1eqn2}
\frac{((1-\alpha)-b(1-\alpha\lambda))\,\Gamma(b)}{\Gamma(1-|a|+b)}\leq\frac{((1-\alpha)-c(1-\alpha\lambda)))\,\Gamma(c)}{\Gamma(1-|a|+c)}
\eeq
\ethm
\begin{proof}  Let $f(z)=z\, _3F_2\left(^{a,b,c}_{b+1,c+1};z\right)$, then, by Lemma \ref{chp2lem1eqn1}, it is enough to show that
\begin{eqnarray*}
\mathcal{T}_1(\alpha,\lambda) &=& \sum_{n=2}^{\infty} [n-(1+n\lambda-\lambda)\alpha]\,|A_n| \leq \alpha-1
\end{eqnarray*}
Using the fact $|(a)_n|\leq  (|a|)_n$, one can get
\begin{eqnarray*}
\mathcal{T}_1(\alpha,\lambda)&=& \sum_{n=2}^{\infty}\,[n(1-\alpha\lambda)-\alpha(1-\lambda)]\,\left(\frac{(|a|)_{n-1}\left(b\right)_{n-1}\left(c\right)_{n-1}}
{\left(b+1\right)_{n-1}\left(c+1\right)_{n-1} \left(1\right)_{n-1}}\right)\\ \\
&=& (1-\alpha\lambda)\, \sum_{n=2}^{\infty} \,n\, \left(\frac{(|a|)_{n-1}\left(b\right)_{n-1}\left(c\right)_{n-1}}
{\left(b+1\right)_{n-1}\left(c+1\right)_{n-1} \left(1\right)_{n-1}}\right)\\
&&\qquad-\alpha\,(1-\lambda)\, \sum_{n=2}^{\infty} \,\left(\frac{(|a|)_{n-1}\left(b\right)_{n-1}\left(c\right)_{n-1}}
{\left(b+1\right)_{n-1}\left(c+1\right)_{n-1} \left(1\right)_{n-1}}\right)\\ \\
&=& (1-\alpha\lambda)\, \sum_{n=0}^{\infty} \,\left(\frac{(n+1)\,(|a|)_{n}\left(b\right)_{n}\left(c\right)_{n}}
{\left(b+1\right)_{n}\left(c+1\right)_{n} \left(1\right)_{n}}\right)-(1-\alpha\lambda)\\
&&\qquad-\alpha\,(1-\lambda)\, \sum_{n=0}^{\infty} \, \left(\frac{(|a|)_{n}\left(b\right)_{n}\left(c\right)_{n}}
{\left(b+1\right)_{n}\left(c+1\right)_{n} \left(1\right)_{n}}\right)+\alpha\,(1-\lambda)
\end{eqnarray*}
Using  the result (1) of Lemma \ref{chp2lem3eqn1} and the formula (\ref{chp2inteq6}) in above mentioned equation, we derived that
\begin{eqnarray*}
&=& (1-\alpha\lambda)\, \frac{bc\, \Gamma(1-|a|) }{c-b}\left[\frac{(1-b)\Gamma(b)}{\Gamma(1-|a|+b)}-\frac{(1-c)\Gamma(c)}{\Gamma(1-|a|+c)}\right]\\
&&\qquad-\alpha\,(1-\lambda)\, \frac{bc\Gamma(1-|a|)}{c-b}\left[\frac{\Gamma(b)}{\Gamma(1-|a|+b)}-\frac{\Gamma(c)}{\Gamma(1-|a|+c)}\right]+\alpha-1\\ \\
&=&  \frac{bc\, \Gamma(1-|a|) }{c-b}\bigg[\frac{(1-b)(1-\alpha\lambda)\,\Gamma(b)}{\Gamma(1-|a|+b)}-\frac{(1-c)(1-\alpha\lambda)\,\Gamma(c)}{\Gamma(1-|a|+c)}\cr
&&\qquad\qquad\qquad\qquad-\frac{\alpha\,(1-\lambda)\,\Gamma(b)}{\Gamma(1-|a|+b)}+\frac{\alpha\,(1-\lambda)\,\Gamma(c)}{\Gamma(1-|a|+c)}\bigg]+\alpha-1\\ \\
&=&  \frac{bc\, \Gamma(1-|a|) }{c-b}\bigg[\frac{(1-\alpha)-b(1-\alpha\lambda))\,\Gamma(b)}{\Gamma(1-|a|+b)}-\frac{((1-\alpha)-c(1-\alpha\lambda))\,\Gamma(c)}{\Gamma(1-|a|+c)}\bigg]+\alpha-1\\
\end{eqnarray*}
The above expression is bounded above by $\alpha-1$ if and only if the equation (\ref{chp2thm1eqn2}) holds, which completes proof.
\end{proof}
\bthm \label{chp2thm2eqn1}
Let $a \in {\IC} \backslash \{ 0 \}, \, b,\,c > 0$,\, $c\neq b$  and $ |a| < \min\{1,\, b+1, \, c+1\}$. A sufficient condition for the function $z\, _3F_2\left(^{a,b,c}_{b+1,c+1};z\right) $ to belong to the class $\N^{\ast}(\lambda,\alpha),\, 1<\alpha \leq \frac{4}{3}$ and $ 0 \leq \lambda  < 1 $ is that
\begin{eqnarray}\label{chp2thm2eqn2}
\frac{(b-1)(b(1-\alpha\lambda)-(1-\alpha))\Gamma(b)}{\Gamma(1-|a|+b)}&\leq&\frac{(c-1)\, (c(1-\alpha\lambda)-(1-\alpha))\Gamma(c)}{\Gamma(1-|a|+c)}
\end{eqnarray}
\ethm
\begin{proof}  Let $f(z)=z\, _3F_2\left(^{a,b,c}_{b+1,c+1};z\right)$, then, by the Lemma \ref{chp2lem2eqn1}, it is enough to show that
\begin{eqnarray*}
\mathcal{T}_2(\alpha,\lambda)  &=& \sum_{n=2}^{\infty} \, n\, [n-(1+n\lambda-\lambda)\alpha]\,|A_n| \leq \alpha-1
\end{eqnarray*}
Using the fact $|(a)_n|\leq  (|a|)_n$, one can get
\begin{eqnarray*}
\mathcal{T}_2(\alpha,\lambda)&=& \sum_{n=2}^{\infty}\, n\, [n(1-\alpha\lambda)-\alpha(1-\lambda)]\,\left(\frac{(|a|)_{n-1}\left(b\right)_{n-1}\left(c\right)_{n-1}}
{\left(b+1\right)_{n-1}\left(c+1\right)_{n-1} \left(1\right)_{n-1}}\right)\\ \\
&=& \sum_{n=2}^{\infty} [n^2\, (1-\alpha\lambda)-\alpha(1-\lambda)\, n]\,\left(\frac{(|a|)_{n-1}\left(b\right)_{n-1}\left(c\right)_{n-1}}
{\left(b+1\right)_{n-1}\left(c+1\right)_{n-1} \left(1\right)_{n-1}}\right)
\end{eqnarray*}
Replace $n = (n-1)+1$ and $n^2  = (n-1)(n-2)+3(n-1)+1$ in above, we find that
\begin{eqnarray*}
\mathcal{T}_2(\alpha,\lambda)&=& \sum_{n=2}^{\infty} [((n-1)(n-2)+3(n-1)+1)]\left(\frac{ (1-\alpha\lambda)\,(|a|)_{n-1}\left(b\right)_{n-1}\left(c\right)_{n-1}}
{\left(b+1\right)_{n-1}\left(c+1\right)_{n-1} \left(1\right)_{n-1}}\right)\\
&&\qquad\qquad-\sum_{n=2}^{\infty} [\alpha(1-\lambda)\, ((n-1)+1)]\left(\frac{(|a|)_{n-1}\left(b\right)_{n-1}\left(c\right)_{n-1}}
{\left(b+1\right)_{n-1}\left(c+1\right)_{n-1} \left(1\right)_{n-1}}\right)\\ \\
&=& (1-\alpha\lambda)\, \sum_{n=2}^{\infty} \, \left(\frac{(n-1)(n-2)\,(|a|)_{n-1}\left(b\right)_{n-1}\left(c\right)_{n-1}}
{\left(b+1\right)_{n-1}\left(c+1\right)_{n-1} \left(1\right)_{n-1}}\right)\\
&&\qquad\qquad+(3-2\alpha\,\lambda-\alpha)\, \sum_{n=2}^{\infty} \, \left(\frac{(n-1)\,(|a|)_{n-1}\left(b\right)_{n-1}\left(c\right)_{n-1}}
{\left(b+1\right)_{n-1}\left(c+1\right)_{n-1} \left(1\right)_{n-1}}\right)\\
&&\qquad\qquad\qquad\qquad\qquad+(1-\alpha)\, \sum_{n=2}^{\infty} \,\left(\frac{(|a|)_{n-1}\left(b\right)_{n-1}\left(c\right)_{n-1}}
{\left(b+1\right)_{n-1}\left(c+1\right)_{n-1} \left(1\right)_{n-1}}\right)\\ \\
&=& (1-\alpha\lambda)\, \sum_{n=3}^{\infty} \, \left(\frac{\,(|a|)_{n-1}\left(b\right)_{n-1}\left(c\right)_{n-1}}
{\left(b+1\right)_{n-1}\left(c+1\right)_{n-1} \left(1\right)_{n-3}}\right)\\
&&\qquad\qquad+(3-2\alpha\,\lambda-\alpha)\, \sum_{n=2}^{\infty} \, \left(\frac{(|a|)_{n-1}\left(b\right)_{n-1}\left(c\right)_{n-1}}
{\left(b+1\right)_{n-1}\left(c+1\right)_{n-1} \left(1\right)_{n-2}}\right)\\
&&\qquad\qquad\qquad\qquad\qquad+(1-\alpha)\, \sum_{n=2}^{\infty} \,\left(\frac{(|a|)_{n-1}\left(b\right)_{n-1}\left(c\right)_{n-1}}
{\left(b+1\right)_{n-1}\left(c+1\right)_{n-1} \left(1\right)_{n-1}}\right)\\ \\
&=& (1-\alpha\lambda)\, \sum_{n=0}^{\infty} \, \left(\frac{\,(|a|)_{n+2}\left(b\right)_{n+2}\left(c\right)_{n+2}}
{\left(b+1\right)_{n+2}\left(c+1\right)_{n+2} \left(1\right)_{n}}\right)\\
&&\qquad\qquad+(3-2\alpha\,\lambda-\alpha)\, \sum_{n=0}^{\infty} \, \left(\frac{(|a|)_{n+1}\left(b\right)_{n+1}\left(c\right)_{n+1}}
{\left(b+1\right)_{n+1}\left(c+1\right)_{n+1} \left(1\right)_{n}}\right)\\
&&\qquad\qquad\qquad\qquad\qquad+(1-\alpha)\, \sum_{n=1}^{\infty} \,\left(\frac{(|a|)_{n}\left(b\right)_{n}\left(c\right)_{n}}
{\left(b+1\right)_{n}\left(c+1\right)_{n} \left(1\right)_{n}}\right)\\ \\
&=& (1-\alpha\lambda)\left(\frac{|a|(|a|+1)b(b+1)c(c+1)}{(b+1)(b+2)(c+1)(c+2)}\right) \sum_{n=0}^{\infty} \, \left(\frac{(|a|+2)_{n}\left(b+2\right)_{n}\left(c+2\right)_{n}}
{\left(b+3\right)_{n}\left(c+3\right)_{n} \left(1\right)_{n}}\right)\\
&&\qquad+(3-2\alpha\,\lambda-\alpha)\, \left(\frac{abc}{(b+1)(c+1)}\right)\, \sum_{n=0}^{\infty} \, \left(\frac{(|a|)_{n+1}\left(b\right)_{n+1}\left(c\right)_{n+1}}
{\left(b+1\right)_{n+1}\left(c+1\right)_{n+1} \left(1\right)_{n}}\right)\\
&&\qquad\qquad\qquad\qquad+(1-\alpha)\, \sum_{n=0}^{\infty} \,\left(\frac{(|a|)_{n}\left(b\right)_{n}\left(c\right)_{n}}
{\left(b+1\right)_{n}\left(c+1\right)_{n} \left(1\right)_{n}}\right)-(1-\alpha)
\end{eqnarray*}
Using the formula (\ref{chp2inteq6}) in above mentioned equation, we find that
\begin{eqnarray*}
&=& (1-\alpha\lambda)\,\left(\frac{|a|(|a|+1)b(b+1)c(c+1)}{(b+1)(b+2)(c+1)(c+2)}\right)\,  \left(\frac{(b+2)(c+2)\Gamma{(1-(a+2))}}{(c+2)-(b+2)}\right)\\
&&\qquad\times \left(\frac{\Gamma(b+2)}{1-(|a|+2)+(b+2)}-\frac{\Gamma(c+2)}{1-(|a|+2)+(c+2)}\right)\\
&&+(3-2\alpha\,\lambda-\alpha)\, \left(\frac{|a|bc}{(b+1)(c+1)}\right)\, \left(\frac{(b+1)(c+1)\Gamma{(1-(|a|+1))}}{(c+1)-(b+1)}\right)\\
&&\qquad\times \left(\frac{\Gamma(b+1)}{1-(|a|+1)+(b+1)}-\frac{\Gamma(c+1)}{1-(|a|+1)+(c+1)}\right)\\
&&+(1-\alpha)\, \frac{bc\Gamma(1-|a|)}{c-b}\left[\frac{\Gamma(b)}{\Gamma(1-|a|+b)}-\frac{\Gamma(c)}{\Gamma(1-|a|+c)}\right]-(1-\alpha)\\ \\
&=& (1-\alpha\lambda)\, \left(\frac{bc\,(-|a|)(-(|a|+1))\,\Gamma{(1-(|a|+2))}}{c-b}\right) \left(\frac{(b+1)\,b\,\Gamma(b)}{1-|a|+b}-\frac{(c+1)\,c\,\Gamma(c)}{1-|a|+c}\right)\\
&&-(3-2\alpha\,\lambda-\alpha)\, \left(\frac{bc (-|a|)\Gamma{(1-(|a|+1))}}{c-b}\right)\left(\frac{b\,\Gamma(b)}{1-|a|+b}-\frac{c\,\Gamma(c)}{1-|a|+c}\right)\\
&&+(1-\alpha)\, \frac{bc\Gamma(1-|a|)}{c-b}\left[\frac{\Gamma(b)}{\Gamma(1-|a|+b)}-\frac{\Gamma(c)}{\Gamma(1-|a|+c)}\right]-(1-\alpha)
\end{eqnarray*}
Using $\Gamma(1-a)=-a\Gamma(-a)$, the aforesaid equation reduces to
\begin{eqnarray*}
&=& \left(\frac{bc\, \Gamma(1-|a|)}{c-b}\right)\\
&&\times\left[\frac{(b-1)(b(1-\alpha\lambda)-(1-\alpha))\Gamma(b)}{\Gamma(1-|a|+b)}-\frac{(c-1)\, (c(1-\alpha\lambda)-(1-\alpha))\Gamma(c)}{\Gamma(1-|a|+c)}\right]+\alpha-1
\end{eqnarray*}
The above expression is bounded above by $\alpha-1$ if and only if the equation (\ref{chp2thm2eqn2}) holds, which completes proof.
\end{proof}
\blem\label{chp2lem4eqn1}\cite{chp2dixitpal1995}
If $f\in \R^{\tau}(A,B)$ is of the form (\ref{chp2inteq0}), then
\begin{eqnarray}\label{chp2lem4eqn2}
  |a_n| &\leq& (A-B)\frac{|\tau|}{n},\, n\in \IN \smallsetminus \{1\}.
\end{eqnarray} The result is sharp.
\elem
Using the Lemma \ref{chp2lem4eqn1}, we prove the following results:
\bthm\label{chp2thm3eqn0}
Let $a \in {\IC} \backslash \{ 0 \}, \, b,\,c > 0$,\, $c\neq b$  and $ |a| < \min\{1,\, b+1, \, c+1\}$ and $f\in \R^{\tau}(A,B)\cap \V$. Then $\mathcal{I}_{b+1, c+1}^{a,b,c}(f)(z) \in \mathcal{N}^{\ast}(\alpha, \lambda)$ if
\begin{eqnarray} \label{chp2thm3eqn1}
\bigg( \frac{bc\, \Gamma(1-|a|) }{c-b}\bigg[\frac{(1-\alpha)-b(1-\alpha\lambda))\,\Gamma(b)}{\Gamma(1-|a|+b)}-\frac{((1-\alpha)-c(1-\alpha\lambda))\,\Gamma(c)}{\Gamma(1-|a|+c)}\bigg]\bigg)\nonumber&&\\
\qquad \times\left(\frac{(A-B)\,|\tau|}{(1-(A-B)\,|\tau|)}\,\right)&\leq& \alpha-1.
\end{eqnarray}
\ethm
\begin{proof}  Let $f$ be of the form (\ref{chp2inteq0}) belong to the class $\R^{\tau}(A,B)\cap \V$. Because of Lemma \ref{chp2lem2eqn1}, it is enough to show that
\begin{eqnarray*}
&& \sum_{n=2}^{\infty}\, n\, [n(1-\alpha\lambda)-\alpha(1-\lambda)]\,\left(\frac{(|a|)_{n-1}\left(b\right)_{n-1}\left(c\right)_{n-1}}
{\left(b+1\right)_{n-1}\left(c+1\right)_{n-1} \left(1\right)_{n-1}}\right)|a_n| \leq \alpha-1
\end{eqnarray*}
since $f\in \R^{\tau}(A,B)\cap \V$, then by Lemma \ref{chp2lem4eqn1}, we have $$|a_n| \leq (A-B)\frac{|\tau|}{n},\, n\in \IN \smallsetminus \{1\}. $$
Letting
\begin{eqnarray*}
\mathcal{T}_3(\alpha,\lambda)  &=& \sum_{n=2}^{\infty}\, n\, [n(1-\alpha\lambda)-\alpha(1-\lambda)]\,\left(\frac{(|a|)_{n-1}\left(b\right)_{n-1}\left(c\right)_{n-1}}
{\left(b+1\right)_{n-1}\left(c+1\right)_{n-1} \left(1\right)_{n-1}}\right)|a_n|
\end{eqnarray*}
we derived that
\begin{eqnarray*}
\mathcal{T}_3(\alpha,\lambda)  &=& (A-B)\,|\tau|\,\sum_{n=2}^{\infty}\, [n(1-\alpha\lambda)-\alpha(1-\lambda)]\,\left(\frac{(|a|)_{n-1}\left(b\right)_{n-1}\left(c\right)_{n-1}}
{\left(b+1\right)_{n-1}\left(c+1\right)_{n-1} \left(1\right)_{n-1}}\right)\\ \\
&=& (A-B)\,|\tau|\,\bigg((1-\alpha\lambda)\, \sum_{n=2}^{\infty} \,n\, \left(\frac{(|a|)_{n-1}\left(b\right)_{n-1}\left(c\right)_{n-1}}
{\left(b+1\right)_{n-1}\left(c+1\right)_{n-1} \left(1\right)_{n-1}}\right)\cr
&&\qquad-\alpha\,(1-\lambda)\, \sum_{n=2}^{\infty} \,\left(\frac{(|a|)_{n-1}\left(b\right)_{n-1}\left(c\right)_{n-1}}
{\left(b+1\right)_{n-1}\left(c+1\right)_{n-1} \left(1\right)_{n-1}}\right)\bigg)\\ \\
&=& (A-B)\,|\tau|\,\bigg((1-\alpha\lambda)\, \sum_{n=0}^{\infty} \,\left(\frac{(n+1)\,(|a|)_{n}\left(b\right)_{n}\left(c\right)_{n}}
{\left(b+1\right)_{n}\left(c+1\right)_{n} \left(1\right)_{n}}\right)-(1-\alpha\lambda)\cr
&&\qquad-\alpha\,(1-\lambda)\, \sum_{n=0}^{\infty} \, \left(\frac{(|a|)_{n}\left(b\right)_{n}\left(c\right)_{n}}
{\left(b+1\right)_{n}\left(c+1\right)_{n} \left(1\right)_{n}}\right)+\alpha\,(1-\lambda)\bigg)
\end{eqnarray*}
Using the result (1) of Lemma \ref{chp2lem3eqn1} and the formula (\ref{chp2inteq6}) in above mentioned equation, we derived that
\begin{eqnarray*}
&=&(A-B)\,|\tau|\,\bigg( (1-\alpha\lambda)\, \frac{bc\, \Gamma(1-|a|) }{c-b}\left[\frac{(1-b)\Gamma(b)}{\Gamma(1-|a|+b)}-\frac{(1-c)\Gamma(c)}{\Gamma(1-|a|+c)}\right]\\
&&\qquad-\alpha\,(1-\lambda)\, \frac{bc\Gamma(1-|a|)}{c-b}\left[\frac{\Gamma(b)}{\Gamma(1-|a|+b)}-\frac{\Gamma(c)}{\Gamma(1-|a|+c)}\right]+\alpha-1\bigg)\\ \\
&=&(A-B)\,|\tau|\,\bigg(  \frac{bc\, \Gamma(1-|a|) }{c-b}\bigg[\frac{(1-b)(1-\alpha\lambda)\,\Gamma(b)}{\Gamma(1-|a|+b)}-\frac{(1-c)(1-\alpha\lambda)\,\Gamma(c)}{\Gamma(1-|a|+c)}\cr
&&\qquad\qquad\qquad\qquad-\frac{\alpha\,(1-\lambda)\,\Gamma(b)}{\Gamma(1-|a|+b)}+\frac{\alpha\,(1-\lambda)\,\Gamma(c)}{\Gamma(1-|a|+c)}\bigg]+\alpha-1\bigg)\\ \\
&=& (A-B)\,|\tau|\,\bigg( \frac{bc\, \Gamma(1-|a|) }{c-b}\bigg[\frac{(1-\alpha)-b(1-\alpha\lambda))\,\Gamma(b)}{\Gamma(1-|a|+b)}-\frac{((1-\alpha)-c(1-\alpha\lambda))\,\Gamma(c)}{\Gamma(1-|a|+c)}\bigg]\cr &&\qquad+\alpha-1\bigg)
\end{eqnarray*}
The above expression is bounded above by $\alpha-1$ if and only if the equation (\ref{chp2thm3eqn1}) holds, which completes proof.
\end{proof}
\bthm\label{chp2thm4eqn0}
Let $a \in {\IC} \backslash \{ 0 \}, \, b,\,c > 0$,\, $c\neq b$  and $ |a| < \min\{1,\, b+1, \, c+1\}$ and $f\in \R^{\tau}(A,B)\cap \V$. Then $\mathcal{I}_{b+1, c+1}^{a,b,c}(f)(z) \in \mathcal{M}^{\ast}(\alpha, \lambda)$ if
\begin{eqnarray} \label{chp2thm4eqn1}
\bigg(\frac{ (1-\alpha\lambda)\, bc\, \Gamma(1-|a|) }{c-b}\left[\frac{\Gamma(b)}{\Gamma(1-|a|+b)}-\frac{\Gamma(c)}{\Gamma(1-|a|+c)}\right]\qquad\qquad\qquad\qquad\qquad\qquad\qquad&&\cr
-\left(\frac{\alpha\,(1-\lambda)\,bc}{(|a|-1)(b-1)(c-1)}\right)\left[\frac{\Gamma(2-|a|)}{c-b}\left(\frac{(c-1)\Gamma(b)}{\Gamma(1-|a|+b)}-\frac{(b-1)\Gamma(c)}{\Gamma(1-|a|+c)}\right)-1\right]\bigg)\nonumber&&\\
\times\left(\frac{(A-B)\,|\tau|}{(1-(A-B)\,|\tau|)}\,\right)\leq \alpha-1.&&
\end{eqnarray}
\ethm
\begin{proof}  Let $f$ be of the form (\ref{chp2inteq0}) belong to the class $\R^{\tau}(A,B)\cap \V$. Because of Lemma \ref{chp2lem1eqn1}, it is enough to show that
\begin{eqnarray*}
&& \sum_{n=2}^{\infty}\,  [n(1-\alpha\lambda)-\alpha(1-\lambda)]\,\left(\frac{(|a|)_{n-1}\left(b\right)_{n-1}\left(c\right)_{n-1}}
{\left(b+1\right)_{n-1}\left(c+1\right)_{n-1} \left(1\right)_{n-1}}\right)|a_n| \leq \alpha-1
\end{eqnarray*}
since $f\in \R^{\tau}(A,B)\cap \V$, then by Lemma \ref{chp2lem4eqn1} the inequality (\ref{chp2lem4eqn2}) holds. Letting
\begin{eqnarray*}
\mathcal{T}_4(\alpha,\lambda)  &=& \sum_{n=2}^{\infty}\, [n(1-\alpha\lambda)-\alpha(1-\lambda)]\,\left(\frac{(|a|)_{n-1}\left(b\right)_{n-1}\left(c\right)_{n-1}}
{\left(b+1\right)_{n-1}\left(c+1\right)_{n-1} \left(1\right)_{n-1}}\right)|a_n|
\end{eqnarray*}
We get
\begin{eqnarray*}
\mathcal{T}_4(\alpha,\lambda)  &=& (A-B)\,|\tau|\,\sum_{n=2}^{\infty}\, \frac{1}{n} \,[n(1-\alpha\lambda)-\alpha(1-\lambda)]\,\left(\frac{(|a|)_{n-1}\left(b\right)_{n-1}\left(c\right)_{n-1}}
{\left(b+1\right)_{n-1}\left(c+1\right)_{n-1} \left(1\right)_{n-1}}\right)\\ \\
&=& (A-B)\,|\tau|\,\bigg((1-\alpha\lambda)\, \sum_{n=2}^{\infty} \, \left(\frac{(|a|)_{n-1}\left(b\right)_{n-1}\left(c\right)_{n-1}}
{\left(b+1\right)_{n-1}\left(c+1\right)_{n-1} \left(1\right)_{n-1}}\right)\cr
&&\qquad-\alpha\,(1-\lambda)\, \sum_{n=2}^{\infty}\,\frac{1}{n} \,\left(\frac{(|a|)_{n-1}\left(b\right)_{n-1}\left(c\right)_{n-1}}
{\left(b+1\right)_{n-1}\left(c+1\right)_{n-1} \left(1\right)_{n-1}}\right)\bigg)\\ \\
&=& (A-B)\,|\tau|\,\bigg((1-\alpha\lambda)\, \sum_{n=0}^{\infty} \,\left(\frac{(|a|)_{n}\left(b\right)_{n}\left(c\right)_{n}}
{\left(b+1\right)_{n}\left(c+1\right)_{n} \left(1\right)_{n}}\right)-(1-\alpha\lambda)\cr
&&\qquad-\alpha\,(1-\lambda)\, \sum_{n=0}^{\infty} \, \left(\frac{(|a|)_{n}\left(b\right)_{n}\left(c\right)_{n}}
{\left(b+1\right)_{n}\left(c+1\right)_{n} \left(1\right)_{n+1}}\right)+\alpha\,(1-\lambda)\bigg)
\end{eqnarray*}
Using the formula (\ref{chp2inteq6}) and the result (4) of Lemma \ref{chp2lem3eqn1} in above mentioned equation, we have
\begin{eqnarray*}
&=&(A-B)\,|\tau|\,\bigg( (1-\alpha\lambda)\, \frac{bc\, \Gamma(1-|a|) }{c-b}\left[\frac{\Gamma(b)}{\Gamma(1-a+b)}-\frac{\Gamma(c)}{\Gamma(1-|a|+c)}\right]\\
&&\, -\left(\frac{\alpha\,(1-\lambda)\,bc}{(|a|-1)(b-1)(c-1)}\right)\left[\frac{\Gamma(2-|a|)}{c-b}\left(\frac{(c-1)\Gamma(b)}{\Gamma(1-|a|+b)}-\frac{(b-1)\Gamma(c)}{\Gamma(1-|a|+c)}\right)-1\right]\cr&&\qquad\qquad\qquad+\alpha-1\bigg)
\end{eqnarray*}
The above expression is bounded above by $\alpha-1$ if and only if the equation (\ref{chp2thm4eqn1}) holds, which completes proof.
\end{proof}

\end{document}